\documentclass[12pt]{amsart}
\usepackage[latin1]{inputenc}

\newtheorem{thm}{Theorem}

\newtheorem{lem}[thm]{Lemma}

\theoremstyle{definition}
\newtheorem{defn}[thm]{Definition}
\newtheorem{prob}{Problem}

 \numberwithin{equation}{section}

\newcommand{\To}{\longrightarrow}
\begin{document}

\thanks{This research was partially supported by the grant
BFM2002-01719 of MCyT (Spain) and a FPU grant of MEC (Spain).} \subjclass[2000]{46B50, 46B26, 54B10, 54D30}
\keywords{Uniform Eberlein compact, regular averaging operator, countable product}
\title[]{Countable products of spaces of finite sets}%
\author{Antonio Avilés}

\address{Departamento de Matemáticas\\ Universidad de Murcia\\ 30100 Espinardo (Murcia)\\ Spain }%
\email{avileslo@um.es}

\begin{abstract}
We consider the compact spaces $\sigma_n(\Gamma)$ of subsets of
$\Gamma$ of cardinality at most $n$ and their countable products.
We give a complete classification of their Banach spaces of
continuous
functions and a partial topological classification.\\
\end{abstract}

\maketitle

For an infinite set $\Gamma$ and a natural number $n$, we consider
the space
$$\sigma_n(\Gamma) = \{ x\in \{0,1\}^\Gamma : |supp(x)|\leq n\}.$$

Here $supp(x) = \{\gamma\in\Gamma : x_\gamma\neq 0\}$. This is a
closed, hence compact subset of $\{0,1\}^\Gamma$, which is
identified with the family of all subsets of $\Gamma$ of
cardinality at most $n$. In this work we will study the spaces
which are countable products of spaces $\sigma_n(\Gamma)$, mainly
their topological classification as well as the classification of
their Banach spaces of continuous
functions.\\

Let $T$ be the set of all sequences $(\tau_n)_{n=1}^\infty$ with $0\leq\tau_n\leq\omega$. When $\tau$ runs over
$T$, $\sigma_\tau(\Gamma) = \prod_1^\infty\sigma_n(\Gamma)^{\tau_n}$ runs over all finite and countable products
of spaces $\sigma_k(\Gamma)$. For $\tau\in T$ we will call $j(\tau)$ to the supremum of all $n$ with $\tau_n>0$
and $i(\tau)$ to the supremum of all $n$ with $\tau_n = \omega$. If $\tau_n<\omega$ for all $n\geq 1$, then
$i(\tau)=0$. Always $0\leq i(\tau)\leq j(\tau)\leq\omega$. Theorem \ref{clasificacion topologica} below
summarizes our knowledge about the topological classification and its proof consists of a number of lemmas along
Section~\ref{secttopo}.

\begin{thm}\label{clasificacion topologica}
Let $\tau,\tau'\in T$ and $\Gamma$ an uncountable set.
\begin{enumerate} \item Suppose $j(\tau)<\omega$. In this case,
$\sigma_{\tau'}(\Gamma)$ is homeomorphic to $\sigma_\tau(\Gamma)$ if and only if $i(\tau)=i(\tau')$ and
$\tau_n=\tau'_n$ for all $n>i(\tau)$. \item Suppose $i(\tau)=\omega$. In this case, if $i(\tau')=\omega$, then
$\sigma_\tau(\Gamma)$ is homeomorphic to $\sigma_{\tau'}(\Gamma)$.
\end{enumerate}
\end{thm}

This is not a complete classification and leaves the following question open:

\begin{prob}
Let $\Gamma$ be an uncountable set and $\tau,\tau'\in T$ such that $j(\tau')=j(\tau)=\omega$, $i(\tau)<\omega$
and there exists some $n\geq i(\tau)$ with $\tau_n\neq\tau'_n$. Is $\sigma_{\tau}(\Gamma)$ homeomorphic to
$\sigma_{\tau'}(\Gamma)$?
\end{prob}

For example, one particular instance of the problem is whether $\prod_{i=1}^\infty\sigma_i(\Gamma)$ is
homeomorphic to
$\prod_{i=2}^\infty\sigma_i(\Gamma)$.\\

About the spaces of continuous functions, it has been recently proved by Marciszewski \cite{Marciszewskic0} that
a Banach space $C(K)$ is isomorphic to $c_0(\Gamma)$ if and only if $K\subset\sigma_n(\Gamma)$ for some
$n<\omega$. This is the case of any compact of the form $K=\prod_{i=1}^{n}\sigma_{k_i}(\Gamma)$ which can be
embedded into $\sigma_{\sum k_i}(\bigcup_1^n\Gamma\times\{i\})$ by $x\mapsto \bigcup_1^n x_i\times\{i\}$. Hence,
it is a consequence of Marciszewski's result that the Banach spaces of continuous functions over finite products
of spaces $\sigma_k(\Gamma)$ over a fixed $\Gamma$ are all isomorphic. In Section~\ref{sectBanach} we prove a
similar result for countable products:

\begin{thm}\label{classifproducts}
Let $\Gamma$ be an infinite set and $(k_n)$ be any sequence of positive integers. Then the Banach spaces
$C(\prod_{n<\omega}\sigma_{k_n}(\Gamma))$ and $C(\sigma_{1}(\Gamma)^\omega)$ are isomorphic.
\end{thm}

The techniques that we will use are based on the use of regular averaging operators and the so called Pe\l
czy\'nski's decomposition method, developed in \cite{Miljutin} and \cite{Pelczynskiexaves} in order to achieve
Miljutin's result that the spaces of continuous functions over uncountable metrizable compacta are all
isomorphic.

\begin{defn}
 Let $\phi:L\To K$ be a continuous surjection between compact spaces.
 A regular averaging operator for $\phi$ is a bounded positive linear operator $T:C(L)\To
C(K)$ with $T(1_L)=1_K$ and $T(x\circ \phi)=x$ for all $x\in C(K)$.
\end{defn}

The countable products of spaces of the form $\sigma_n(\Gamma)$ are uniform Eberlein compact spaces, cf.
\cite{BenSta}. This class consists of the weakly compact subsets of the Hilbert spaces, or equivalently of the
compact subsets of the space $$B(\Gamma) = \{ x\in [-1,1]^\Gamma : \sum_{\gamma\in\Gamma}|x_\gamma|\leq 1\}\sim
(B_{\ell_2(\Gamma)},w)$$ for some set $\Gamma$. Indeed, $\sigma_n(\Gamma)$ is homeomorphic to
$B(\Gamma)\cap\{0,\frac{1}{n}\}^\Gamma$. We establish the following result:

\begin{thm}\label{RAOBRW}
Let $K$ be a uniform Eberlein compact of weight $\kappa$. Then there is a closed subspace $L$ of
$\sigma_1(\kappa)^\mathbb{N}$ and an onto continuous map $f:L\To K$ which admits a regular averaging operator.
\end{thm}

This improves a result of Argyros and Arvanitakis~\cite{ArgArvRAO} that for every uniform Eberlein compact space
$K$ there is a totally disconnected uniform Eberlein compact space $L$ of the same weight and a continuous
surjection $f:L\To K$ which admits a regular averaging operator, and also a result of Benyamini, Rudin and Wage
\cite{BenRudWag} that every uniform Eberlein compact of weight $\kappa$ is a continuous image of a closed subset
of $\sigma_1(\kappa)^\mathbb{N}$. We note that there are many totally disconnected uniform Eberlein compact
spaces which cannot be embedded into $\sigma_1(\kappa)^\mathbb{N}$, cf. Lemma~\ref{embebimientos de sigmaenes}
below.

\section*{Notations}

All topological spaces will be assumed to be completely regular.
By identifying elements of $\{0,1\}^\Gamma$ with subsets of
$\Gamma$, the space $\sigma_n(\Gamma)\subset\{0,1\}^\Gamma$ can be
viewed as the family of all subsets of $\Gamma$ of cardinality
less than or equal to $n$, endowed with the topology which has a
base the sets of the form

$$\Phi_F^G = \{y\in\sigma_{n}(\Gamma) : F\subset y \subset \Gamma\setminus G\}$$

for $F$ and $G$ finite subsets of $\Gamma$. We will denote by
$p:\sigma_1(\Gamma)^k\To\sigma_k(\Gamma)$ the continuous
surjection given by $p(x_1,\ldots,x_k)=x_1\cup\cdots\cup x_k$.
Note that from the existence of such a function follows the fact
that any countable product $\prod_{i<\omega}\sigma_{k_i}(\Gamma)$
is a continuous image of $\sigma_1(\Gamma)^\omega$. We will also
denote $B^+(\Gamma)=B(\Gamma)\cap [0,1]^\Gamma$.

\section{Banach space classification}\label{sectBanach}

The following Theorem \ref{p admite RAO} is the key result of this section. A somewhat similar fact can be found
in \cite{Shchepinspectra}, namely that the natural surjection $K^2\To exp_2(K)=\{\{x,y\} : x,y\in K\}$ given by
$(x,y)\mapsto \{x,y\}$ has a regular averaging operator.

\begin{thm}\label{p admite RAO}
The map $p:\sigma_1(\Gamma)^k\To\sigma_k(\Gamma)$ admits a regular
averaging operator.
\end{thm}

Proof: For every $y\in\sigma_k(\Gamma)$ let us denote by $L(y)$ the subset of $p^{-1}(y)$ consisting of all
$(x^1,\ldots,x^k)\in p^{-1}(y)$ such that $x^i\cap x^j =\emptyset$ for $i\neq j$ (that is, $L(y)$ consists of
those tuples of $p^{-1}(y)$ in which no singleton appears twice).

 The regular
averaging operator $T:C(\sigma_1(\Gamma)^k)\To
C(\sigma_k(\Gamma))$ is defined as follows:

$$ T(f)(y) = \frac{1}{|L(y)|}\sum_{x\in L(y)}f(x) $$

The only difficult point is in proving that $T(f)$ is a continuous
function whenever $f$ is continuous. So fix $f\in
C(\sigma_1(\Gamma)^k)$ and a point $y\in\sigma_k(\Gamma)$ and
$\varepsilon>0$. For each $x=(x_1,\ldots,x_k)\in L(y)$, since $f$
is continuous at $x$, there is a neighborhood $U_x$ of $x$ in
$\sigma_1(\Gamma)^k$ in which $\sup_{x'\in
U_x}|f(x)-f(x')|<\varepsilon$. The set $U_x$ must contain a basic
neighborhood of $x$ of the form

$$\Phi_{x_1}^{G^x_1}\times\cdots\times\Phi_{x_k}^{G^x_k}\subset U_x$$

where $G^x_i$ is a finite set of $\Gamma$ disjoint with $x_i$. We
define a neighborhood of $y$ as

$$V = \Phi_y^{\bigcup_{x\in L(y)}\bigcup_{i=1}^{k}G_i^x\setminus y}$$

and we shall see that $|T(f)(y)-T(f)(y')|<\varepsilon$ for every $y'\in V$. So we fix $y'\in V$ (in particular
$y\subset y'$). First, we define an onto map $r:L(y')\To L(y)$ in the following way, if $ (x_1,\ldots,x_k)\in
L(y')$ then $r(x)= (r(x)_1,\ldots,r(x)_k)$ where $r(x)_i = x_i\cap y$. It is straightforward to check that all
the fibers of $r$ have the same cardinality, call $n=|r^{-1}(x)|$, so that $|L(y')| = n|L(y)|$. The key fact
(used in the final inequality in the expression below) is that if $x\in L(y)$ and $x'\in r^{-1}(x)$, then $x'\in
U_x$. To see this, take $x=(x_1,\ldots,x_k)\in L(y)$ and $x'=(x'_1,\ldots,x'_k)\in r^{-1}(x)$. We check that
$x'_i\in\Phi_{x_i}^{G_i^x}$. If $x'_i \subset y$ then $x'_i=x_i$. If $x'_i = \{\gamma\}\subset y'\setminus y$
then $x_i=\emptyset$ and since $y'\in V$, $\gamma\not\in G_i^x$ and again $x'_i\in\Phi_{x_i}^{G_i^x}$. Finally,
\begin{eqnarray*}
|T(f)(y') - T(f)(y)| &=& \left|\frac{1}{|L(y')|}\sum_{x'\in
L(y')}f(x') -
\frac{1}{|L(y)|}\sum_{x\in L(y)}f(x)\right|\\
&=& \left|\frac{1}{|L(y')|}\sum_{x\in L(y)}\sum_{x'\in
r^{-1}(x)}f(x') - \frac{1}{|L(y)|}\sum_{x\in L(y)}f(x)\right|\\
&=& \left|\frac{1}{n|L(y)|}\sum_{x\in L(y)}\sum_{x'\in
r^{-1}(x)}f(x') - \frac{1}{|L(y)|}\sum_{x\in L(y)}f(x)\right|\\
&=& \left|\frac{1}{|L(y)|}\sum_{x\in
L(y)}\left(\left(\frac{1}{n}\sum_{x'\in
r^{-1}(x)}f(x')\right) - f(x)\right)\right|\\
&=& \left|\frac{1}{|L(y)|}\sum_{x\in
L(y)}\left(\frac{1}{n}\sum_{x'\in
r^{-1}(x)}(f(x') - f(x))\right)\right|\\
&\leq& \frac{1}{|L(y)|}\sum_{x\in
L(y)}\left(\frac{1}{n}\sum_{x'\in
r^{-1}(x)}|f(x') - f(x)|\right)\\
&<& \frac{1}{|L(y)|}\sum_{x\in L(y)}\left(\frac{1}{n}\sum_{x'\in
r^{-1}(x)}\varepsilon\right) = \varepsilon\\
\end{eqnarray*}
\qed

\begin{lem}\label{propiedades conocidas de RAO}
\begin{itemize}
\item[(a)] Let $g:L\To K$ be a continuous surjection between compact spaces which admits a regular averaging
operator and let $M$ be a closed subset of $K$. Then the restriction $g:g^{-1}(M)\To M$ also admits a regular
averaging operator~\cite[Proposition 18]{ArgArvRAO}. \item[(b)] Let $\{g_i:L_i\To K_i\}$ be a family of
continuous surjections between compact spaces which admit regular averaging operators. Then the product map
$\prod g_i:\prod L_i\To \prod K_i$ admits a regular averaging operator too~\cite[Proposition
4.7]{Pelczynskiexaves}.
\end{itemize}
\end{lem}

Proof of Theorem \ref{RAOBRW}: We make the observation that the space $B(\Gamma)$ can be embedded into
$B^{+}(\Gamma\times\{a,b\})\sim B^{+}(\Gamma)$ by the map $u(x)_{\gamma,a}= \max(0,x_\gamma)$ and
$u(x)_{\gamma,b} = \max(0,-x_\gamma)$. This observation allows to consider our $K$ as a subset of $B^+(\Gamma)$
with $|\Gamma|=\kappa$. Let $\phi:\{0,1\}^\omega\To [0,1]$ given by $\phi(x)=\sum r_i x_i$ where $r_i=
\frac{1}{3}\left(\frac{2}{3}\right)^i$. It is proven in~\cite{ArgArvRAO} that $\phi$ admits a regular averaging
operator and hence by Lemma~\ref{propiedades conocidas de RAO} also $\phi^\Gamma:\{0,1\}^{\omega\times\Gamma}\To
[0,1]^\Gamma$ and its restriction $\phi^\Gamma:L'=(\phi^\Gamma)^{-1}(K)\To K$ admit a regular averaging
operator. The space $L'$ is a subspace of $L_0= (\phi^\Gamma)^{-1}(B^+(\Gamma))$ for which we can give the
following description:
\begin{eqnarray*} x\in L_0 &\iff& \phi^\Gamma(x)\in B^+(\Gamma)\\
&\iff& \sum_{\gamma\in\Gamma}\phi^\Gamma(x)_\gamma \leq 1\\
&\iff& \sum_{\gamma\in\Gamma}\sum_{n=0}^\infty r_n x_{(\gamma,n)}\leq 1\\
&\iff& \sum_{n=0}^\infty r_n N_n(x)\leq 1,\\
\end{eqnarray*}

where $N_n(x)$ is the cardinality of $supp(x|_{\Gamma\times\{n\}})$. From this description, if $M_n$ denotes the
integer part of $r_n^{-1}$, then $L'\subset L_0\subset \prod_{n=1}^\infty\sigma_{M_n}(\Gamma)$. From
Theorem~\ref{p admite RAO} and part (b) of Lemma~\ref{propiedades conocidas de RAO} follows the existence of a
continuous surjection $g:\sigma_1(\Gamma)^\omega\To\prod_{n=1}^\infty\sigma_{M_n}(\Gamma)$ which admits a
regular averaging operator. Making use of part (a) of Lemma~\ref{propiedades conocidas de RAO} we get a
surjection $g:L=g^{-1}(L')\To L'$ with regular averaging operator and the composition $L\To L'\To K$ is the
desired map.$\qed$\\

We shall need now the so called Pe\l czy\'nski's decomposition
method, which is used to establish the existence of isomorphisms
between Banach spaces. For Banach spaces $X$ and $Y$ we shall
write $X|Y$ if there exists a Banach space $Z$ such that $X\oplus
Z$ is isomorphic to $Y$, shortly $X\oplus Z\sim Y$. Also,
$Y=(X_1\oplus X_2\oplus\cdots)_{c_0}$ denotes the $c_0$-sum of the
Banach spaces $X_1,X_2,\ldots$, $$Y=\{ y=(x_n)\in\prod X_n :
\lim\| x_n\|=0\}, \ \ \|y\|=\sup_n\|x_n\|.$$

\begin{thm} [cf. \cite{Pelczynskiexaves}, \S 8]\label{Pelczdecmeth} Let $X$ and $Y$ be Banach spaces such that
$X|Y$, $Y|X$ and $(X\oplus X\oplus\cdots)_{c_0}\sim X$, then
$X\sim Y$.
\end{thm}

If there exists a surjection $\phi:L\To K$ with regular averaging operator, then $C(K)|C(L)$, cf.
\cite{Pelczynskiexaves}. In particular if $L\subset K$ is a retract of $K$, since in this case the restriction
operator is a regular averaging operator for the retraction. On the other hand, in order to guarantee the last
hypothesis in Theorem~\ref{Pelczdecmeth} we shall use the criterion of Lemma~\ref{c0 sumas de compactos}. For
topological spaces $K_n$, $K_1\oplus K_2\oplus\cdots$ denotes the discrete topological sum, while $\alpha(S)$ is
the one point compactification of a locally compact space $S$.

\begin{lem}\label{c0 sumas de compactos}
Let $K$ be a compact space which is homeomorphic to $\alpha(K\oplus K\oplus\cdots)$. Then $(C(K)\oplus
C(K)\oplus\cdots)_{c_0}\sim C(K)$.
\end{lem}

Proof: We apply Theorem~\ref{Pelczdecmeth} to $X=(C(K)\oplus C(K)\oplus\cdots)_{c_0}$ and $Y=C(K)$. The only
point is in checking that $X|Y$. Let $\infty$ denote the infinity point of $\alpha(K\oplus K\oplus\cdots)\sim
K$. Then $X\sim Y'=\{f\in C(K) : f(\infty)=0\}$ and $Y\sim Y'\oplus\mathbb{R}$.$\qed$\\

Proof of Theorem~\ref{classifproducts}:  Set
$K=\sigma_1(\Gamma)^\omega$ and $L=\prod\sigma_{k_n}(\Gamma)$. We
apply Theorem~\ref{Pelczdecmeth} to $X=C(K)$ and $Y = C(L)$.
First, we already observed that from Theorem~\ref{p admite RAO}
and Lemma~\ref{propiedades conocidas de RAO}(b) follows the
existence of a surjection with regular averaging operator $f:K\To
L$ and hence $C(L)|C(K)$. On the other hand, $K$ is a retract of
$L$ because for any $k$, $\sigma_1(\Gamma)$ is homeomorphic to a
clopen subset of $\sigma_{k}(\Gamma)$: the family of all subsets
which contain fixed elements $\gamma_1,\ldots,\gamma_{k-1}$.
Therefore $C(K)|C(L)$. By Lemma~\ref{c0 sumas de compactos}, it
only remains to show that $\alpha(K\oplus K\oplus\cdots)\sim K$.
For this, fix $\gamma\in\Gamma$ and set for $n=1,2,\ldots$

$$K_n = \{x\in K=\sigma_1(\Gamma)^\omega : \gamma\in
x_1\cap\cdots\cap x_{n-1}\setminus x_n\}.$$

The sets $K_n$ are disjoint clopen sets homeomorphic to $K$ and $K$ is the one point compactification of their
union with point of infinity $(\{\gamma\},\{\gamma\},\ldots)$.$\qed$

\section{Topological classification}\label{secttopo}

This section is devoted to the proof of Theorem~\ref{clasificacion topologica}. Before entering this, we point
out why we assume $\Gamma$ to be uncountable. The reasonings below do not apply in the countable case and the
situation is indeed completely different. All perfect totally disconnected metrizable compact spaces are
homeomorphic \cite[Theorem 7.4]{Kechrisbook} and this implies that all countable products of spaces
$\sigma_k(\omega)$ are homeomorphic. The finite products are countable compacta, whose topological
classification is also well known after the classical paper \cite{MazSie}: two of them are homeomorphic if and
only if they have same Cantor-Bendixson derivation index and the same cardinality of the last nonempty
Cantor-Bendixson derivative. Straightforward computations give that these two invariants for a finite product
$\prod_{i=1}^n\sigma_{k_i}(\omega)$ take the values $1+\sum_1^n k_i$ and 1 respectively. From now on, $\Gamma$
will be always an uncountable set.

\begin{lem}\label{absorcion chica}
If $m<n$ then $\sigma_m(\Gamma)\times\sigma_n(\Gamma)^\omega$ is homeomorphic to $\sigma_n(\Gamma)^\omega$.
\end{lem}

Proof: We denote again by $(X_1\oplus X_2\oplus\cdots)$ the discrete topological sum of the spaces $X_1, X_2,
\ldots$ and by $\alpha X$ the one-point compactification of the locally compact space $X$. Fix
$\gamma_0,\ldots,\gamma_{n-1}\in\Gamma$. We consider the set $L=\omega\times\{0,\ldots,n-1\}$ endowed with the
lexicographical order: $(k,i)<(k',i')$ whenever either $k<k'$ or $k=k'$ and $i<i'$. For every $(k,i)\in L$ we
define a clopen set of $\sigma_n(\Gamma)^\omega$ as
\begin{eqnarray*} A_{(k,i)} &=& \{x\in\sigma_n(\Gamma)^\omega : \gamma_i\not\in
x_k, \gamma_{i'}\in x_{k'} \forall (k',i')<(k,i)\}\\
&=& \{x\in\sigma_n(\Gamma)^\omega : \gamma_i\not\in
x_k\supset\{\gamma_0,\ldots,\gamma_{i-1}\} , x_j =
\{\gamma_0,\ldots,\gamma_{n-1}\} \forall j<k\}.
\end{eqnarray*}

Notice that $A_{(k,i)}$ is homeomorphic to $\sigma_{n-i}(\Gamma)\times\sigma_n(\Gamma)^\omega$ and that $\{A_l :
l\in L\}$ constitutes a disjoint sequence of clopen subsets of $\sigma_n(\Gamma)^\omega$ with only limit point
the sequence $\xi\in\sigma_n(\Gamma)^\omega$ constantly equal to $\{\gamma_0,\ldots,\gamma_{n-1}\}$. Hence,
$$\sigma_n(\Gamma)^\omega\approx\alpha\left(\bigoplus_{l\in
L}A_l\right)\approx\alpha\left(\bigoplus_{i=0}^{n-1}\bigoplus_{j<\omega}
\left(\sigma_{n-i}(\Gamma)\times\sigma_n(\Gamma)^\omega\right)\right).$$

On the other hand, we can perform a similar decomposition in $\sigma_m(\Gamma)\times\sigma_n(\Gamma)^\omega$
defining, for $j< m$ and $(k,i)\in L$:
\begin{eqnarray*}
B'_{j} &=& \{(y,x)\in\sigma_m(\Gamma)\times\sigma_n(\Gamma)^\omega
: \gamma_j\not\in y, \{\gamma_0,\ldots,\gamma_{j-1}\}\subset
y\}\\
B_{(k,i)} &=&
\{(y,x)\in\sigma_m(\Gamma)\times\sigma_n(\Gamma)^\omega :
\gamma_i\not\in x_k, \gamma_{i'}\in x_{k'} \forall
(k',i')<(k,i),\\ & & \ \ \ \{\gamma_0,\ldots,\gamma_{m-1}\}\subset
y\}\end{eqnarray*}

Again $B'_j$ is homeomorphic to $\sigma_{m-j}(\Gamma)\times\sigma_n(\Gamma)^\omega$, $B_{(k,i)}$ is homeomorphic
to $\sigma_{n-i}(\Gamma)\times\sigma_n(\Gamma)^\omega$ and altogether they constitute a disjoint sequence of
clopen sets with a single limit point $(\{\gamma_0,\ldots,\gamma_{m-1}\},\xi)$ out of them, so
$$\sigma_m(\Gamma)\times\sigma_n(\Gamma)^\omega\approx\alpha\left(\bigoplus_{l\in
L}B_l\oplus\bigoplus_{j=0}^{m-1}B'_j\right)\approx\alpha\left(\bigoplus_{i=0}^{n-1}\bigoplus_{j<\omega}
\left(\sigma_{n-i}(\Gamma)\times\sigma_n(\Gamma)^\omega\right)\right).$$

\begin{lem}\label{absorcion grande}
If $m<n<\omega$ then $\sigma_m(\Gamma)^\omega\times\sigma_n(\Gamma)^\omega$ is homeomorphic to
$\sigma_n(\Gamma)^\omega$.
\end{lem}

Proof: $\sigma_m(\Gamma)^\omega\times\sigma_n(\Gamma)^\omega\approx
\left(\sigma_m(\Gamma)\times\sigma_n(\Gamma)^\omega\right)^\omega \approx
\left(\sigma_n(\Gamma)^\omega\right)^\omega\approx \sigma_n(\Gamma)^\omega$. $\qed$

\begin{lem}\label{absorcion general}
Let $m_1,\ldots,m_r<n<\omega$ and $e_1,\ldots,e_r\leq\omega$. Then
the space $\prod_{i=1}^r\sigma_{m_i}(\Gamma)^{e_i}\times
\sigma_n(\Gamma)^\omega$ is homeomorphic to
$\sigma_n(\Gamma)^\omega$.
\end{lem}

Proof: Follows from repeated application of Lemmas~\ref{absorcion
chica} and~\ref{absorcion grande} above.$\qed$\\

From Lemma~\ref{absorcion general} it follows that any space $\sigma_\tau(\Gamma)$ with $i(\tau)=\omega$ is
homeomorphic to $\sigma_{(\omega,\omega,\ldots)}(\Gamma)$ (because we can substitute each factor
$\sigma_n(\Gamma)^\omega$ of $\sigma_\tau(\Gamma)$ by the homeomorphic $\prod_{i\leq n}\sigma_i(\Gamma)^\omega$)
and this proves part (2) of Theorem~\ref{clasificacion topologica}. Lemma~\ref{absorcion general} also shows
that it is irrelevant in determining the homeomorphism class of $\sigma_\tau(\Gamma)$ which are the values
$\tau_n$ for $n<i(\tau)$. Hence, in order to prove part(1) of Theorem~\ref{clasificacion topologica} it remains
to show that if $j(\tau)<\omega$ and $\sigma_\tau(\Gamma)$ is homeomorphic to $\sigma_{\tau'}(\Gamma)$ then
$\tau_n=\tau'_n$ for all $n>i(\tau)$.\\

We recall that a family $\{S_\eta\}_{\eta\in H}$ of sets is a
$\Delta$-system if there is a set $S$ (called the root of the
$\Delta$-system) such that $S_\eta\cap S_{\eta'} = S$ for all
$\eta\neq\eta'$. We will make use of the fact that any uncountable
family of finite sets has an uncountable subfamily which is a
$\Delta$-system, cf.\cite[Theorem 1.4]{Chainconditions} for
$\kappa=\omega$ and $\alpha=\omega_1$.\\

The following lemma includes as a particular case that $\sigma_{n+1}(\Lambda)$ does not embed into
$\sigma_n(\Gamma)^\omega$. This fact, whose proof corresponds to Steps 1 - 3 below was shown to us by Witold
Marciszewski, and it seems that it was known to several people before.

\begin{lem}\label{embebimientos de sigmaenes}
If $|\Lambda|>\omega$, $n\geq 0$, $k\geq 0$, then the space $\sigma_{n+1}(\Lambda)^{k+1}$ does not embed into
$\sigma_n(\Gamma)^\omega\times\sigma_{n+1}(\Gamma)^k$.
\end{lem}

Proof: Suppose that there exists such an embedding.\\

\emph{Step 1}. Passing to a suitable uncountable subset of
$\Lambda$, we can suppose that there is an embedding
$$\phi:\sigma_{n+1}(\Lambda)^{k+1}\To\sigma_n(\Gamma)^m\times\sigma_{n+1}(\Gamma)^k$$
for some $m<\omega$. To see this, take
$\varphi:\sigma_{n+1}(\Lambda)^{k+1}\To\sigma_n(\Gamma)^\omega\times\sigma_{n+1}(\Gamma)^k$ our original
embedding. In this step, we shall denote an element $x\in\sigma_{n+1}(\Lambda)^{k+1}$ as $x=(x_0,\ldots,x_k)$.
For each $\lambda\in\Lambda$ and every $i\in\{0,\ldots,k\}$ we find a clopen set $A^i_\lambda$ of
$\sigma_n(\Gamma)^\omega\times\sigma_{n+1}(\Gamma)^k$ which separates the disjoint compact sets $\varphi(\{x :
\lambda\in x_i\})$ and $\varphi(\{x : \lambda\not\in x_i\})$. Associated to $A^i_\lambda$ we have a finite
subset $F^i_\lambda\subset\omega$ such that $A^i_\lambda = \sigma_n(\Gamma)^{\omega\setminus F^i_\lambda}\times
B^i_\lambda$ with $B^i_\lambda$ a clopen subset of
$\sigma_n(\Gamma)^{F^i_\lambda}\times\sigma_{n+1}(\Gamma)^{k}$. We choose $\Lambda'$ to be an uncountable subset
of $\Lambda$ such that $\bigcup_{i=0}^k F^i_\lambda=\bigcup_{i=0}^k F^i_{\lambda'}=F$ for all
$\lambda,\lambda'\in\Lambda'$ and in this case the composition
$$\sigma_{n+1}(\Lambda')^{k+1}\hookrightarrow \sigma_{n+1}(\Lambda)^{k+1}\To
\sigma_n(\Gamma)^\omega\times\sigma_{n+1}(\Gamma)^k
\To\sigma_n(\Gamma)^F\times\sigma_{n+1}(\Gamma)^k$$

is one-to-one. The reason is that if $x,y\in \sigma_{n+1}(\Lambda')^{k+1}$ are different then there exists
$i\in\{0,\ldots,k\}$ and $\lambda\in\Lambda'$ such that $\lambda\in x_i$
but $\lambda\not\in y_i$ (or viceversa). Then $\phi(x)\in A_\lambda^i$ and $\phi(y)\not\in A_\lambda^i$ so either
the coordinate of $\sigma_{n+1}(\Gamma)^k$ or some coordinate of $F_i^\lambda\subset F$ must be
different for $\phi(x)$ and $\phi(y)$.\\

\emph{Step 2}. For $i=0,\ldots,k$ and $\lambda\in\Lambda$ we define $e^\lambda_i\in\sigma_{n+1}(\Lambda)^{k+1}$
to be the element which has $\{\lambda\}$ in coordinate $i$ and $\emptyset$ in all other coordinates. Each
$\phi(e_i^\lambda)$ will be of the form
$$\phi(e^\lambda_i) = (x^\lambda_i[1],\ldots,x^\lambda_i[m],x^\lambda_i[m+1],\ldots,x^\lambda_i[m+k])$$
 with $x^\lambda_i[j]\in\sigma_n(\Gamma)$ if $j\leq m$ and
 $x^\lambda_i[j]\in\sigma_{n+1}(\Gamma)$ if $m<j\leq m+k$. Passing to a suitable uncountable subset of
 $\Lambda$, we can assume that for every fixed $i\in\{0,\ldots,k\}$ and $j\in\{1,\ldots,m+k\}$ the
 family $\{x^\lambda_i[j]:\lambda\in\Lambda\}$ is a
 $\Delta$-system of root $R_i[j]$ formed by sets of the same cardinality $c_i[j]$.\\

\emph{Step 3}. We claim that for $i=0,\ldots,n$ and $j=1,\ldots,m$, the $\Delta$-system $\{x^\lambda_i[j] :
\lambda\in\Lambda\}$ is constant. Suppose the contrary for some fixed $i\leq n$ and $j\leq m$. Then
$x^\lambda_i[j] = R\cup S^\lambda\in\sigma_{n}(\Gamma)$ where $R\cap S^\lambda=\emptyset$,
$S^\lambda\neq\emptyset$, and $ S^\lambda\cap S^{\lambda'}=\emptyset$ for $\lambda\neq\lambda'$. We consider the
sets
$$A_\lambda =
\{y=(y[1],\ldots,y[m+k])\in\sigma_n(\Gamma)^m\times\sigma_{n+1}(\Gamma)^k : y[j]\supset S^\lambda\}.$$ The
$A_\lambda$'s are neighborhoods of the $\phi(e^\lambda_i)$'s with the property that for every $F\subset\Lambda$
with $|F|>n$, $\bigcap_{\lambda\in F}A_\lambda = \emptyset$ (because for $y$ in that intersection, $|y[j]|> n$
and $y[j]\in\sigma_n(\Gamma)$). Let $\psi:\sigma_{n+1}(\Lambda)\To\sigma_{n+1}(\Lambda)^{k+1}$ be the map
defined by $\psi(x)_i=x$ and $\psi(x)_{i'}(x)=\emptyset$ if $i'\neq i$. Then the $(\phi\psi)^{-1}(A_\lambda)$'s
are neighborhoods of the $\{\lambda\}$'s in $\sigma_{n+1}(\Lambda)$ with the property that for every
$F\subset\Lambda$ with $|F|>n$, $\bigcap_{\lambda\in F}(\phi\psi)^{-1}(A_\lambda) = \emptyset$. This is a
contradiction since such a family of neighborhoods cannot be found. Namely, take basic neighborhoods with
$\{\lambda\} \in \Phi_{\{\lambda\}}^{G_\lambda}\subset (\phi\psi)^{-1}(A_\lambda)$ and take
$\Lambda'\subset\Lambda$ uncountable with $\{G_\lambda : \lambda\in\Lambda'\}$ a $\Delta$-system of root $R'$.
Then construct inductively a finite sequence $F=\{\lambda_1,\ldots,\lambda_{n+1}\}\subset\Lambda'\setminus R'$
such that $\lambda_p\not\in\bigcup_{q<p}G_{\lambda_q}$ and
$G_{\lambda_p}\cap\{\lambda_1,\ldots,\lambda_{p-1}\}=\emptyset$ (notice that it is possible to choose such a
$\lambda_p$ because $\{\lambda_1,\ldots,\lambda_{p-1}\}\cap R' = \emptyset$ and hence there are only finitely
many $G_\lambda$'s with $\lambda\in\Lambda'$ and
$G_{\lambda}\cap\{\lambda_1,\ldots,\lambda_{p-1}\}\neq\emptyset$). In this case we have
$F\in\bigcap_{\lambda\in F}(\phi\psi)^{-1}(A_\lambda)$.\\

\emph{Step 4}. Notice that, in the case when $k=0$ we already
arrived to a contradiction and the proof is complete. When $k>0$
we need some extra work. From step 3, we deduce that for each
$i\in\{0,\ldots,k\}$ there must exist $j\in\{m+1,\ldots,m+k\}$
such that the family $\{x^\lambda_i[j]: \lambda\in\Lambda\}$ is a
nonconstant $\Delta$-system. Since $i$ runs in a set of $k+1$
elements and $j$ in a set of $k$ elements, there must exist two
different $i,i'\in\{0,\ldots,k\}$ such that for the same $j$,
$\{x^\lambda_i[j]: \lambda\in\Lambda\}$ and $\{x^\lambda_{i'}[j]:
\lambda\in\Lambda\}$ are nonconstant $\Delta$-systems. We assume
that $c_i[j]\geq c_{i'}[j]$ (these numbers are defined in step 2).
 Again, for
$\lambda\in\Lambda$ we consider the sets
$$ A_\lambda = \{(y[1],\ldots,y[m+k])\in\sigma_n(\Gamma)^m\times\sigma_{n+1}(\Gamma)^k
: y[j]\supset x_i^\lambda[j]\},$$
$$ A'_\lambda = \{(y[1],\ldots,y[m+k])\in\sigma_n(\Gamma)^m\times\sigma_{n+1}(\Gamma)^k
: y[j]\supset x_{i'}^\lambda[j]\}.$$

The $A_\lambda$'s and the $A'_\lambda$'s are neighborhoods of the
$\phi(e_i^\lambda)$'s and the $\phi(e_{i'}^\lambda)$'s
respectively with the property that
$$ (\ast)\ \forall\lambda\in\Lambda\forall F\subset\Lambda \left(|F|>n
\wedge x_i^\lambda[j]\not\subseteq\bigcup_{\mu\in F}x_{i'}^\mu[j]\right)
\Rightarrow A_\lambda\cap\bigcap_{\mu\in F}A'_\mu=\emptyset.$$

That intersection is empty because if $y$ belongs to it, then
$$x_i^\lambda[j]\cup\bigcup_{\mu\in F}x_{i'}^\mu[j]\subset
y[j]\in\sigma_{n+1}(\Gamma)$$ and the set in the left, if $x_i^\lambda[j]\not\subseteq\bigcup_{\mu\in
F}x_{i'}^\mu[j]$,  has cardinality greater than $n+1$, a contradiction. Since the $\Delta$-systems are not
constant and $c_i[j]\geq c_{i'}[j]$, if $x_i^\lambda[j]\subseteq\bigcup_{\mu\in F}x_{i'}^\mu[j]$ holds, there
must be some $\mu\in F$ and some $\gamma\in x_i^\lambda[j]$ such that $\gamma\in x_{i'}^\mu[j]\setminus
R_{i'}[j]$. For a fixed $\lambda$ there are only finitely many $\mu$'s with $(x_{i'}^\mu[j]\setminus
R_{i'}[j])\cap x^\lambda_i[j]\neq\emptyset$. Hence for every $\lambda$, we can find a cofinite subset
$\Lambda_\lambda$ of $\Lambda$ such that the hypothesis $x_i^\lambda[j]\not\subseteq\bigcup_{\mu\in
F}x_{i'}^\mu[j]$ of statement $(\ast)$ holds whenever $F\subset\Lambda_\lambda$. For short, we know that for
every $\lambda\in\Lambda$ there exists a cofinite subset $\Lambda_\lambda$ of $\Lambda$ such that
$$\forall F\subset\Lambda_\lambda\ \ |F|>n
\Rightarrow A_\lambda\cap\bigcap_{\mu\in F}A'_\mu=\emptyset.$$ This contradicts the following lemma for
$B_\lambda = \phi^{-1}(A_\lambda)$ and $B'_\lambda = \phi^{-1}(A'_\lambda)$:

\begin{lem}
For every $\lambda\in\Lambda$, let $B_\lambda$ and $B'_\lambda$ be
neighborhoods of $e^\lambda_i$ and $e^\lambda_{i'}$ respectively
in $\sigma_{n+1}(\Lambda)^{k+1}$. Then there exists
$\lambda_0\in\Lambda$ and an infinite set $S\subset\Lambda$ such
that for every $F\subset S$ with $|F|=n+1$,
$$B_{\lambda_0}\cap\bigcap_{\mu\in F}B'_\mu\neq\emptyset$$
\end{lem}

Proof: For a simpler notation, we will assume that $i=0$ and
$i'=1$. Notice that a basic clopen set $\Phi_F^G$ of
$\sigma_{n+1}(\Lambda)$ is nonempty if and only if $F\cap
G=\emptyset$ and $|F|\leq n+1$. Each $B_\lambda$ and each $B'_\mu$
contain basic clopen sets of the form
$$ \Phi^{G^\lambda_0}_{\{\lambda\}}\times
\Phi^{G^\lambda_1}_\emptyset\times\Phi^{G^\lambda_2}_\emptyset\times\cdots\times \Phi^{G^\lambda_k}_\emptyset
\subseteq B_\lambda$$
$$ \Phi^{H^\mu_0}_\emptyset\times
\Phi^{H^\mu_1}_{\{\mu\}}\times\Phi^{H^\mu_2}_\emptyset\times\cdots\times
\Phi^{H^\mu_k}_\emptyset \subseteq B'_\mu$$

with all $G^\lambda_l$ and $H^\mu_l$ finite subsets of $\Lambda$ and $\lambda\not\in G^\lambda_0$ and
$\mu\not\in H^\mu_1$. First, we find $M\subset\Lambda$ a countably infinite set such that $\mu'\not\in H_1^\mu$
for every $\mu,\mu'\in M$. This can be done as follows. We begin with an infinite $M_1\subset\Lambda$ such that
the family $\{H_1^\mu :\mu\in M_1\}$ is a $\Delta$-system of root $R$, and we set $M_2 = M_1\setminus R$. Then
we can find recursively a sequence $(\mu_p)_{p<\omega}\subset M_2$ such that
$\mu_p\not\in\bigcup_{q<p}H_1^{\mu_q}$ and $H_1^{\mu_p}\cap\{\mu_1,\ldots,\mu_{p-1}\}=\emptyset$. After this, we
set $M=\{\mu_p : p<\omega\}$. Now, we choose $\lambda_0\not\in\bigcup_{\mu\in M}H_0^\mu$. Taking $S=\{\mu\in M :
\mu\not\in G_1^{\lambda_0}\}$, then $\lambda_0$ and $S$ are as desired. Namely, take $F\subset S$ with
$|F|=n+1$, and for every $j=0,\ldots,k$ call $I_j = G_j^{\lambda_0}\bigcup_{\mu\in F}H_j^\mu$ so that
$$B_{\lambda_0}\cap\bigcap_{\mu\in F}B'_\mu \supset
\Phi_{\{\lambda_0\}}^{I_0^\mu}\times \Phi_{F}^{I_1^\mu}\times
\prod_{j=2}^{k}\Phi_{\emptyset}^{I_j^\mu}.$$ On the one hand,
$\Phi_{\{\lambda_0\}}^{I_0^\mu}\neq\emptyset$ because we chose
$\lambda_0\not\in\bigcup_{\mu\in M}H_0^\mu$, so $\lambda_0\not\in
I_0^\mu$. On the other hand, $\Phi_{F}^{I_1^\mu}\neq\emptyset$
because, first, since $F\subset M$ and $\mu'\not\in H_1^\mu$ for
every $\mu,\mu'\in M$, it follows that $F\cap\bigcup_{\mu\in
F}H_1^\mu = \emptyset$ and second, since $F\subset S$, just by the
definition of $S$, $F\cap G_1^{\lambda_0} = \emptyset$.$\qed$\\

It is a consequence of Lemma~\ref{embebimientos de sigmaenes} that
$j(\tau)=j(\tau')$ whenever
$\sigma_\tau(\Gamma)=\sigma_{\tau'}(\Gamma)$, since it shows that
$j(\tau)=\omega$ if and only if $\sigma_n(\Gamma)$ can be embedded
into $\sigma_\tau(\Gamma)$ for all $n<\omega$ and, if it is not
the case, $j(\tau)$ is the greatest integer $n$ for which
$\sigma_n(\Gamma)$ embeds into $\sigma_\tau(\Gamma)$. Hence, in
the situation of part (1) of Theorem~\ref{clasificacion
topologica}, it happens that $j(\tau)=j(\tau')=j$ and moreover
 that $\tau_n=\tau'_n$
for all $n\geq j$ since, by Lemma~\ref{embebimientos de sigmaenes} again, $\tau_{j}=\tau'_{j}$ is the greatest
integer $k$ such that $\sigma_j(\Gamma)^k$ embeds into $\sigma_\tau(\Gamma)$ and of course, $\tau_n=\tau'_n=0$
for all $n>j$. In order to finish the proof of this part (1), we must check that $i(\tau)=i(\tau')=i$ and that
$\tau_k = \tau'_k$ for $i<k<j$. In order to get this, we shall look at embeddability of spaces
$\sigma_n(\Gamma)^k$ into the clopen sets of $\sigma_\tau(\Gamma)$. For this purpose, we observe that it is
enough to look at some basic family of clopen sets, if the others are union of them:

\begin{lem}\label{embebimiento en uniones}
Let $X$ be a compact space and $C_1,\ldots,C_t$ open subsets of $X$. If $\sigma_n(\Lambda)^k$ embeds into
$\bigcup_1^t C_i$, then there exists $i\leq t$ such that $\sigma_n(\Lambda)^k$ embeds into $C_i$.
\end{lem}

Proof: It reduces to prove that whenever we express
$\sigma_n(\Lambda)^k$ as a union of open sets as
$$\sigma_n(\Lambda)^k = C_1\cup\cdots\cup C_t$$ then some $C_i$ must
contain a copy of $\sigma_n(\Lambda)^k$. Pick $i\in\{1,\ldots,t\}$
such that $x_0=(\emptyset,\ldots,\emptyset)\in C_i$. There are
finite sets $G^1,\ldots,G^k$ of $\Lambda$ such that
$$x_0\in \Phi_\emptyset^{G^1}\times\cdots\times\Phi_\emptyset^{G^k}\subset C_i$$

This finishes the proof because
$\Phi_\emptyset^{G^1}\times\cdots\times\Phi_\emptyset^{G^k}$
is homeomorphic to $\sigma_n(\Lambda)^k$.$\qed$\\

Let us denote now by $K=\prod_{s\in S}\sigma_{n_s}(\Gamma)$ any finite or countable product of spaces of type
$\sigma_n(\Gamma)$. Any clopen set of $K$ is a finite union of basic clopen sets of the form
$$C = \prod_{s\in A}\Phi_{F_s}^{G_s}\times\prod_{s\not\in A}\sigma_{n_s}(\Gamma)$$

where $A$ is a finite subset of $S$ and $\Phi_{F_s}^{G_s}$ a basic clopen set of $\sigma_{n_s}(\Gamma)$. Such a
basic clopen set is homeomorphic to
$$(\star)\ C\sim \prod_{s\in A}\sigma_{n_s-|F_s|}(\Gamma)\times\prod_{s\not\in
A}\sigma_{n_s}(\Gamma)$$ Now, after Lemma~\ref{embebimientos de sigmaenes}, Lemma~\ref{embebimiento en uniones}
and the topological description $(\star)$ of the basic clopen sets given above, we are in a position to state
that, in the situation of part (1) of Theorem~\ref{clasificacion topologica}, the following hold:

\begin{itemize}
 \item [(A)] $i(\tau)=i(\tau')=i$ is the greatest integer $n$ such that
 $\sigma_n(\Gamma)$ embeds into any clopen set of
 $\sigma_\tau(\Gamma)$.

 \item [(B)] For $n=j,j-1,j-2,\ldots,i+1$, $\tau_n=\tau'_n$ is the greatest integer
 $k$ such that there is a clopen set $C$ of $\sigma_\tau(\Gamma)$ in which
 $\sigma_{n+1}(\Gamma)$ cannot be embedded, but in which nevertheless
 $\sigma_n(\Gamma)^{k+\sum_{r>n}\tau_r}$ does embed.\\
\end{itemize}

This finishes the proof of Theorem \ref{clasificacion topologica}.
For statement (A), since $\sigma_{i(\tau)}(\Gamma)^\omega$ is one
of the factors of $\sigma_\tau(\Gamma)$, it is clear that still
$\sigma_{i(\tau)}(\Gamma)^\omega$ is a factor of any clopen set
like in $(\star)$. On the other hand, there are only finitely many
factors of type $\sigma_m(\Gamma)$, $m>i(\tau)$ in
$\sigma_\tau(\Gamma)$, hence a clopen set like in $(\star)$ can be
obtained so that all factors in $\prod_{s\in
A}\sigma_{n_s-|F_s|}(\Gamma)\times\prod_{s\not\in
A}\sigma_{n_s}(\Gamma)$ are of the form $\sigma_m(\Gamma)$ with
$m\leq i(\tau)$. By Lemma~\ref{embebimientos de sigmaenes},
$\sigma_k(\Gamma)$ does not embed in such $C$ if $k>i(\tau)$.\\

Statement (B) is proved by ``downwards induction'' starting in $j$ and finishing with $i+1$. We know, by
Lemma~\ref{absorcion general}, that
$$\sigma_\tau(\Gamma) \sim
\sigma_{i}(\Gamma)^\omega\times\prod_{m=i+1}^j\sigma_{m}(\Gamma)^{\tau_m}$$ Now statement $(B)$ for $n=j$ is a
direct consequence of Lemma \ref{embebimientos de sigmaenes} since no clopen can contain $\sigma_{j+1}(\Gamma)$
and the maximal exponent of $\sigma_{j}(\Gamma)$ inside $\sigma_\tau(\Gamma)$ is $\tau_j$. We pass to the case
when $i<n<j$. The ``biggest'' possible basic clopen set $C$ of $\sigma_\tau(\Gamma)$ not containing
$\sigma_{n+1}(\Gamma)$ is obtained by reducing as necessary the factors $\sigma_m(\Gamma)$ with $m>n$:
$$C\sim
\sigma_{i}(\Gamma)^\omega\times\prod_{m=i+1}^{n}\sigma_{m}(\Gamma)^{\tau_m}
\times\prod_{m=n+1}^{j}\sigma_{n}(\Gamma)^{\tau_m}$$ The maximal exponent of $\sigma_{n}(\Gamma)$ in such a $C$
is
$\sum_{m=n}^{j}\sigma_{\tau_{m}}$.$\qed$\\

The present work was written during a visit to the University of
Warsaw. The author wishes to thank their hospitality, specially to
Witold Marciszewski and Roman Pol, and to Rafa\l\ Górak, from
the Polish Academy of Sciences. This work owes very much to the
discussion with them and their suggestions.

\end{document}